\documentclass[leqno]{amsart}
\usepackage{amsmath}        
\usepackage{amsfonts}
\usepackage{enumerate}
\usepackage{amsthm}
\usepackage{amssymb}
\usepackage{amscd}
\usepackage[all]{xy}
\usepackage[titletoc,toc,title]{appendix}
\newtheorem{thm}{Theorem}[section] 
\newtheorem{prop}[thm]{Proposition}
\newtheorem{cor}[thm]{Corollary}
\newtheorem{lemma}[thm]{Lemma}
\newtheorem{addendum}[thm]{Addendum}
                                   
\theoremstyle{definition}

\theoremstyle{definition}

\newtheorem{Rem}[thm]{Remark}

\newtheorem{example}[thm]{Example}

\newtheorem{question}[thm]{Question}
\def\sk1{\vskip 10pt}
\def\newline{\hfil\break}

\def\R{\mathbb R}

\def\del{\partial}

\def\sk{\vskip}

\numberwithin{equation}{section}
\title{On semistability of $CAT(0)$ groups}
\author{Ross Geoghegan and Eric Swenson}
\address{\noindent Ross Geoghegan, Department of Mathematical Sciences,
Binghamton University (SUNY),
Binghamton, NY 13902-6000, USA
\newline
\vskip 3pt
Eric Swenson, Department of Mathematics,
Brigham Young University,
Provo, UT 84602, USA}
\thanks{The second-named author was partially supported by a grant from the Simons Foundation (209403)}
\email{ross@math.binghamton.edu,
eric@math.byu.edu}  

\date{June 15, 2018}
\begin{document}
\fontsize{12}{13pt} \selectfont  
\maketitle
\begin{abstract}
Does every one-ended $CAT(0)$ group have semistable fundamental group at
infinity? As we write, this is an open question.  Let $G$ be such a group acting
geometrically on the proper $CAT(0)$ space $X$. In this paper we show
that in order to establish a positive answer to the question it is
only necessary to check that any two geodesic rays in $X$ are properly
homotopic. We then show that if the answer to the question is negative, with $(G,X)$
a counter-example, then the boundary of $X$, $\del X$ with the cone topology, must
have a weak cut point separating two open subsets. This is of interest because a theorem of Papasoglu
and the second-named author \cite{PS} has established that there cannot be
an example of $(G,X)$ where $\del X$ has a cut point. Thus, the search for a negative answer
comes down to the difference between cut points and weak cut points. We also show that the Tits
ball of radius $\frac{\pi}{2}$ about that weak cut point is a ``cut set" in the sense that it separates
$\del X$. Finally, we observe that if a negative example 
$(G, X)$ exists then $G$ is rank 1.
\end{abstract}
\section{Semistability at infinity}\label{stab}
Let $X$ be a proper $CAT(0)$ space having one end, and let $\del X$ denote its compactifying
boundary. A point of $\del X$ is, by definition, an equivalence class of geodesic rays in $X$ any
two of which are boundedly close\footnote{See \cite {BrHa99} for $CAT(0)$ matters.}. A
simple model for $\del X$ is obtained by choosing a base point $b\in X$: then the space of all
geodesic rays starting at $b$, with the compact-open topology, is homeomorphic to $\del X$.
Since $X$ has one end $\del X$ is connected. Following classical terminology we refer to a
connected compact Hausdorff space as a {\it continuum}; all continua in this paper will be
metrizable.
\vskip 5pt
One says that $X$ has {\it semistable fundamental group at infinity} (or {\it has one strong end}) 
if any two proper rays in $X$ are properly homotopic. We recall that the gap between the words
and the definition is explained as follows: An inverse sequence of groups $\{H_{n},
f^{m}_{n}\}$ is said to be {\it semistable} (or {\it Mittag-Leffler}) if, for each $n$, the images
of the bonding homomorphisms $f^{m}_{n}:H_{m}\to H_{n}$ are the same for all
but finitely many values of $m>n$. The relevance to our $X$ is this: Let $\{K_{n}\}$
be an exhausting sequence of compact subsets of $X$ such that, for all
$n$, $K_{n}$ is a subset of the interior of $K_{n+1}$. Then $\{U_{n}:=X-K_{n}\}$ is a basic
system of neighborhoods of infinity. Choosing a suitably parametrized geodesic base ray
$\omega $ in $X$, we consider the inverse sequence of groups $\{\pi _{1}(U_{n}, \omega),
\text{(inclusion)}_{\#}\}$. {\it It is a fact that this sequence of groups is semistable if and only if
any two proper rays in $X$ are properly homotopic;} see \cite{G}, Section 16.1 for details. 
\vskip 5pt
Extending our simple model of $\del X$ as the space of geodesic
rays beginning at $b$, we define $\widehat X$ to be the space of all
maps $[0,\infty )\to X$ which are {\it either} geodesic rays starting at
$b$, {\it or} are geodesic segments on a finite interval $[0, \rho]$,
mapping $0$ to $b$, and are constant on $[\rho ,\infty)$. Understood
here is the compact-open topology, so, by the Arzela-Ascoli Theorem,
$\widehat X$ is compact. It is only a minor abuse of language to
say that $X$ is open and dense in $\widehat X$. It is obvious that
$\del X$ is a $Z$-set\footnote{This means that for any open set $U$
in $\widehat X$, the inclusion map $U\cap X\to U$ is a homotopy
equivalence. Equivalently, for any positive $\epsilon $, any map of
a polyhedron $P$ into $\widehat X$ can be $\epsilon $-homotoped off
$\del X$, holding fixed an arbitrarily large closed subset of $P$
lying in the pre-image of $X$.} in $\widehat X$. We write ${\widehat
U_{n}}:=U_{n}\cup \del X$. Because $\del X$ is a $Z$-set we obtain
an inverse sequence of open subsets of $\widehat X$ with the following
property:
\vskip 5pt
{\it The ``inclusion map" $\{\pi _{1}(U_{n}, \omega),
\text{(inclusion)}_{\#}\}\to \{\pi _{1}({\widehat U_{n}}, \omega (\infty )),
\text{(inclusion)}_{\#}\}$ induces an isomorphism of pro-groups.}
\vskip 5pt
Thus, semistability at infinity of $X$ can be read off from either inverse sequence. 
\vskip 5pt
We will be dealing with the relationship between $X$ having semistable fundamental group at
infinity and shape properties of $\del X$. There is a delicate issue here which makes necessary a
brief discussion of shape theory. 
\section{Remarks on shape theory}
For our purposes, shape theory of continua can be presented in the following limited way.
Suppose we are given a compact absolute retract (AR) $\widehat W$ and a subcontinuum $Y$ which is
a $Z$-set in $\widehat W$. We choose a base point $y\in Y$. The dense open set 
$W:={\widehat W}-Y$ is a locally compact AR having one end. If $\{V_{n}\}$ is a basic
system of neighborhoods of infinity in $W$ then ${\widehat V_{n}}:=V_{n}\cup Y$ is a basic
system of neighborhoods of $Y$ in $\widehat W$. Moreover:
\begin{enumerate}[(1)]
\item If each $V_n$ is closed in $W$ then  $\{{\widehat V_{n}}\}$ is a basic system of closed
neighborhoods of $Y$ in $\widehat W$;
\item If each $V_n$ is an absolute neighborhood retract (ANR), then $\{{\widehat V_{n}}\}$ is
a basic system of ANR neighborhoods of $Y$ in $\widehat W$ (because ``adding" a $Z$-set
does not destroy the ANR property). 
\end{enumerate}
\vskip 5pt
One really wants both $(1)$ and $(2)$ to hold because then one has a
basic sequence of compact ANR neighborhoods of $Y$ in the AR $\widehat
W$, in which case  any pro-homotopy\footnote{Because of our concern with
fundamental group we are really interested in pointed shape invariance
--- i.e. invariance under pro-homotopy where all maps and homotopies
respect base points. It is cumbersome to keep inserting the word
``pointed", and we shall often omit it.} invariant of the inverse sequence
$\{({\widehat V_{n}},y)\}$ is an intrinsic shape theoretic invariant of
$(Y,y)$. In particular, this means that if $\{{\widehat V'_{n}}\}$
is a sequence of compact ANR neighborhoods of a $Z$-set copy of $Y$
in a (possibly different) AR $W'$, then the same (pointed) pro-homotopy
properties will hold there. One such property is semistability of the
inverse sequence of fundamental groups. If any such inverse sequence of
pointed compact ANR's has semistable fundamental pro-group then $(Y,y)$
is said to have {\it semistable fundamental pro-group}\footnote{Or, in
the language of shape theory, $Y$ is {\it pointed $1$-movable}.}. It is
a fact that (because $Y$ is connected) this property is independent of
the choice of base point.
\vskip 5pt
{\bf Example:} It is well-known that every proper $CAT(0)$ space is an
AR. Since the passage from $X$ to $\widehat X$ is achieved by
``adding" a $Z$-set, namely $\del X$, the compact space $\widehat
X$ is also an AR. Thus $({\widehat X}, \del X)$ is an example of
$({\widehat W},Y)$ of the previous paragraphs. We apply the above discussion to
the neighborhoods $U_n$ of Section \ref{stab}. In order to equate
``$X$ having semistable fundamental group at infinity"  with ``$\del X$
having semistable fundamental pro-group" we would like to be able to
choose the $U_n$ (or a cofinal sequence of neighborhoods $U'_n$) to be
closed ANR subsets of $X$; i.e. to satisfy $(1)$ and $(2)$. When $X$ is triangulable (for
example, a $CAT(0)$ cubical complex) this is easy to achieve. But we do not know how to find
such closed ANR neighborhoods of infinity in an arbitrary proper $CAT(0)$ space $X$.
Nevertheless, we assert:
\begin{prop}\label{tech} A proper $CAT(0)$ space $X$ has semistable fundamental group at
infinity if and only if its boundary $\del X$ has semistable fundamental pro-group.
\end{prop}
\vskip 5pt
Because the proof of Proposition \ref{tech} for the non-triangulable case involves ideas outside
the scope of this paper we defer it to the Appendix.
\section{Geodesic rays}\label{rays}
Since geodesic rays are proper, one may ask if geodesic rays alone determine whether or not $X$
has semistable fundamental group at infinity.
\begin{thm}\label{straight} For a one-ended proper $CAT(0)$ space $X$ the following are
equivalent:
\begin{enumerate}[(i)]
\item $X$ has semistable fundamental group at infinity;
\item any two geodesic rays in $X$ are properly homotopic;
\item $\del X$ has semistable fundamental pro-group.
\end{enumerate}
\end{thm} 
\begin{Rem}\label{straight1} Any two geodesic rays in $X$ are properly
homotopic through geodesic rays (i.e. every level of the homotopy a
geodesic ray) if and only if $\del X$ is path connected. Theorem \ref{straight}
describes something weaker than path connectedness.  \end{Rem}
\begin{Rem} It has long been known (see Kraszinkiewicz \cite {K}) that if a continuum is
locally connected then it has semistable fundamental pro-group, so Theorem \ref{straight} is
mainly of interest when $\del X$ is not known to be locally connected.  
\end{Rem}
\begin{Rem}\label{straight3} That theorem of Krasinkiewicz \cite {K} says more: A continuum
has semistable fundamental pro-group if and only if it is shape equivalent to a locally connected
continuum. 
\end{Rem}
\vskip 5pt
We need some shape theoretic terminology. 
\begin{enumerate}[(1)]
\item A {\it strong shape component} of $\del X$ is a proper homotopy class of proper rays in
$X$. 
\vskip 5pt
\item The proper ray $c:[0,\infty )\to X$ {\it ends at the point} $p\in
\del X$ if the map $c$ extends continuously to ${\hat c}:[0,\infty ]\to
\widehat X$ by mapping the point $\infty $ to $p$. 
\vskip 5pt
\item Two points $p$ and $q$ of $\del X$ lie in the same {\it component of
joinability} if there are proper rays in $X$ ending at $p$ and at $q$
which are properly homotopic.
\end{enumerate} 
\begin{Rem} For strong shape theory see, for example, Section 17.7 of \cite{G}. Components of
joinability were introduced in \cite{KM}.
\end{Rem}
\begin{lemma}\label{LC1} If two proper rays end at the same point of $\del X $ then they are
properly homotopic.
\end{lemma}
\begin{proof} Let the proper rays $c$ and $c'$ end at the point
$p$. AR's are locally simply connected, so we can choose a basic system
of neighborhoods of $p$ in $\widehat X$
\begin{equation*}
\dots \subseteq W_n\subseteq V_n \subseteq U_n=W_{n-1}\subseteq V_{n-1}\subseteq
U_{n-1}\subseteq\dots
\end{equation*}
such that 
\begin{enumerate}[(i)] 
\item $[c(n),\infty )\cup [c'(n), \infty
)\subseteq W_n$; 
\item any two points in $W_n$ can be joined by a path
in $V_n$; 
\item any loop in $V_n$ is homotopically trivial in $U_n$.
\end{enumerate} 
Because $\del X$ is a $Z$-set, any two points in $W_{n}\cap X$ can be joined by a path in $V_{n}\cap X$, and any loop in $V_{n}\cap X$ is homotopically trivial in $U_{n}\cap X$.
For each $n$ choose a path $\omega _n$ in $V_{n}\cap X$ joining
$c(n)$ to $c'(n)$. Choose a trivializing homotopy in $U_{n-1}\cap X$ of the
loop formed by the segments $[c'(n-1), c'(n)]$,  $[c(n-1), c(n)]$,
$\omega _{n-1}$ and $\omega _n$. Together, these give the desired  proper homotopy
between $c$ and $c'$.
\end{proof}
We denote by $[c]$ the strong shape component of $\del X$ defined by
the proper ray $c$. If some $c'$ in $[c]$ ends at a point of $\del X$
we say that the strong shape component $[c]$ of $\del X$ is {\it non-empty}. In
general there may also be empty strong shape components, meaning examples
of $[c]$ containing no such $c'$. The existence of these is explored
in \cite{GK} where, in particular the following is proved (\cite{GK}
Corollary 3.6 and Theorem 5.1):
\begin{thm}\label{kras} If $\del X$ has an empty strong shape component then it has
uncountably many such, and also uncountably many components of joinability.
\end{thm}
\begin{Rem} The proof in \cite{GK} is for metrizable continua in
general, not specifically for boundaries of $CAT(0)$ spaces.  
\end{Rem}
{\it Proof of Theorem \ref{straight}:} The equivalence of $(i)$ and
$(iii)$ is Proposition \ref{tech}. The direction``$(i)$ implies $(ii)$"
is trivial. Assume that any two geodesic rays in $X$ are properly
homotopic. Then $\del X$ has only one component of joinability. By
Theorem \ref{kras}, $\del X$ has no empty strong shape components. So
any strong shape component $[c]$ contains a proper ray $c'$ which
ends at a point $p$ of $\del X$. There is also a geodesic ray ending
at $p$, and, by Lemma \ref{LC1}, it is properly homotopic to  $c'$. It
follows that 
$X$ has semistable fundamental group at infinity.  \hfill$\square$
\section{Relevance to group theory}
A finitely presented one-ended group $G$ {\it has semistable fundamental
group at infinity} if the universal cover of some (equivalently,
any) finite $2$-complex whose fundamental group is isomorphic to $G$
has semistable fundamental group at infinity in the sense defined above
for $X$.   This property of finitely presented groups is a quasi-isometry
invariant; see, for example, Section 18.2 of \cite {G}. It is unknown
if every finitely presented group has this property. All one-ended
hyperbolic groups have it, as do many other classes of groups.
\vskip 5pt
When the discrete group $G$ acts geometrically (i.e. properly discontinuously and
cocompactly) on the one-ended proper $CAT(0)$ space $X$ then $G$ is quasi-isometric to $X$,
so $G$ has one end and  is finitely presented. Moreover, $G$ has semistable fundamental
group at infinity if and only if the same is true of $X$. The relevance to group
theory is that in order to show that a $CAT(0)$ group $G$ has semistable fundamental group at
infinity, one need only check the condition on geodesic rays given in
Theorem \ref{straight}. And it need only be checked on one proper $CAT(0)$ space on which
$G$ acts geometrically.
\vskip 5pt
This is of interest because of the following still-unsolved problem:
\vskip 5pt
\begin{question}\label{mainq}Is it true that every one-ended $CAT(0)$ group has semistable
fundamental group at infinity?
\end{question}
An immediate corollary of Theorem \ref{straight} is: 
\begin{cor}
If $G$ is not rank $1$ then $G$ has semistable fundamental group at infinity.
\end{cor}
\begin{proof} If $G$ is not rank $1$ then, by a theorem of Ballmann
and Buyalo \cite{BaBu}, the Tits diameter of $\del X$ is finite,
so $\del X$ is Tits path connected. The obvious function from the Tits
boundary to $\del X$ is continuous so $\del X$ is path connected. Remark
\ref{straight1} applies.  \end{proof}
The theorem of Krasinkiewicz \cite {K} mentioned in Remark \ref{straight3} implies
that if there exists a $CAT(0)$ group $G$ which does not have semistable fundamental group at
infinity then $\del X$ is not shape equivalent to a locally connected continuum.
\vskip 5pt
\section{Cut points and semistability}
We continue to explore properties which $\del X$ would have to possess if its
fundamental
pro-group were not semistable. 
\vskip 5pt
A point $c$ of a continuum $Y$ is a {\it cut point} if $Y-\{c\}$ is not connected. The
point $c$
is a {\it weak cut point} if there are two other points $a$ and $b$ of $Y$ such that any
subcontinuum containing $a$ and $b$ must also contain $c$. In that case we say that $c$
{\it
weakly separates} $a$ from $b$. A variant on this will be important: let $U$ and $V$ be
disjoint subsets of $Y-\{c\}$; the point $c$ {\it weakly separates} $U$ {\it from} $V$ if
any subcontinuum of $Y$ which meets $U$ and $V$ must also contain $c$.
\begin{example} The limit segment $I$ of the Topologist's Sine Curve $S$ provides an
example: $c$ is the mid point of $I$, $a$ and $b$ are its two end points. This $c$ is a
weak cut point in $S$ but it is not a cut point. However, we note that $c$ does not weakly
separate open neighborhoods of $a$ and of $b$.
\end{example}
Let $X$ be as before, and let $G$ act geometrically on $X$. It is proved in \cite{PS} that
the
existence of such a group $G$ implies that $\del X$ cannot contain a cut point. In this
section we show that if $\del X$ does not have semistable fundamental pro-group then
there are open subsets $U$ and $V$ of $\del X$ and a point $c\in \del X$ such that $c$
weakly separates $U$ from $V$.
\vskip 5pt
A proper $CAT(0)$ space $X$ is {\it almost geodesically complete} if there
is a number $\lambda $ such that for any points $p$ and $q$ of $X$ there
is a geodesic ray starting at $p$ and passing within distance $\lambda$
of $q$. It is proved in \cite{On05} that when (as here) there is a group $G$
acting geometrically on $X$ by isometries, then $X$ is almost geodesically complete. {\it
From
now on this number $\lambda $ will be part of the data coming with $X$.}
\vskip 5pt
\begin{lemma}\label{semi} Let $p\in X$. Suppose that for every $\rho >0$ there
exists $\sigma >\rho $ with the following property: Whenever $\tau \geq \sigma $ and
$\alpha ,
\beta :[0,\infty)\to X$ are unit speed geodesic rays starting at $p$ with
$d(\alpha (\tau), \beta (\tau))<2\lambda +1$ then the geodesic segment $[\alpha (\tau ),
\beta (\tau )]$ can be homotoped, relative to $\alpha$ and $\beta $, out of the closed ball
$B(p,
\tau +\lambda +1)$ by a homotopy missing $B(p,\rho )$. Then $X$ has semistable
fundamental
group at infinity.
\end{lemma}         
\begin{proof}
Let $\omega $ be a geodesic base ray starting at $p$, and let $\rho
>0$ be given.  Choose $\sigma $ as in the hypothesis, and let $E$ be
any compact set containing $D:=B(p,\sigma +2\lambda +1)$. Let $\gamma
:[0,1]\to X$ be a loop in $X- D$ based at $\omega $ and let $\delta
:=d(p,\gamma(0))$.  By very small homotopies, we may assume that $\gamma$
is piecewise geodesic.  Choose a partition $0=x_0,x_1,\dots x_n =1$
with the property that the length of $\gamma \mid :[x_{i-1},x_i] \to X$
is at most $1$.  For each $1\leq i <n$ choose a unit speed geodesic ray
$\omega _{i}$ starting at $p$ with $d(\omega _{i}(y_i), \gamma(x_i))\leq
\lambda $ for some $y_i$ (necessarily, $y_i\geq \sigma +\lambda +1$), setting
$\omega _{0}=\omega _{n} =\omega $, and $y_{0}=\delta =y_n$.  
\vskip 5pt
The loop formed from $\gamma ([x_{i-1},x_i])$ and the geodesic segments
$[\gamma(x_i),\omega _i(y_i)]$, $[\gamma(x_{i-1}),\omega _{i-1}(y_{i-1})]$
and $[\omega _{i-1}(y_{i-1}), \omega _i(y_i)]$ has length bounded by
$4\lambda +2$.   The straight line homotopy of this loop to the point
$\gamma(x_i)$ misses $B(p,\rho)$.  Thus we may homotope $\gamma $
to the piecewise geodesic loop with vertices $\omega _0(y_0), \omega
_1(y_1), \dots ,\omega _{n-1}(y_{n-1})$.  Each segment has length at most
$2\lambda +1$, and $y_i\geq \sigma +\lambda +1$.
This new loop is homotopic to the piecewise geodesic segment with 
vertices $\omega _0(\sigma +\lambda +1), \omega _1(\sigma +\lambda +1),
\dots \omega _{n-1}(\sigma +\lambda +1)$; again, the homotopy misses
$B(p,\rho)$.  By hypothesis, we can homotope this, relative to $\omega
_0,\dots , \omega _{n-1}$,  to a loop $\eta $ supported in $X-B(p,
\sigma +2\lambda +2)$, again missing $B(p,\rho)$.  Now consider the
subpath $\eta_i$  of $\eta$ from $\omega _i$ to $\omega _{i+1}$ (mod
$n$)   Repeating the above process we can homotope each path $\eta_i$
outside of $B(p,\sigma +2\lambda  + 3)$ relative to $\omega _i$ and $\omega
_{i+1}$ missing $B(p, \rho)$.  Gluing these homotopies we get a homotopy of $\gamma
$ to a loop outside $B(p,\sigma +2\lambda + 3)$, where the homotopy is relative to
$\omega $ and is supported outside $B(p,\rho )$. Iterating the process, for any $n$ we can
homotope $\gamma$ out of $B(p,\sigma +2\lambda  +n)$ by a homotopy relative to
$\omega $ missing $B(p,\rho)$.  
Thus we can move $\gamma$ off the compact set $E$.  
\end{proof}
\vskip 5pt
For $Z\subseteq \del X$ and $x\in X$ the {\it cone} with base $Z$ and
vertex $x$ is the subset ${\mathcal C}_{x}Z$ of $\widehat X$ consisting
of $Z$ and the images in $X$ of all the geodesic rays $[0,\infty )\to X$
starting at $x$ and ending in $Z$.
\begin{lemma}\label{compact} When $Z$ is compact and $c\in \del X - Z$ is
defined by the geodesic ray $\gamma$ starting at $x$ then, given $\kappa >0$,
there is $\sigma >0$ such that for all $\tau >\sigma $ the distance in $X$
from $\gamma (\tau)$ to  ${\mathcal C}_{x}Z$  is greater than $\kappa $.
\end{lemma}
\begin{Rem} In fact it can be proved that the cone is compact when $Z$
is compact, but we will not need this.  
\end{Rem}
\begin{proof} ({\it of Lemma \ref{compact}}) Suppose not. Then there
is $\kappa >0$, a sequence $(\xi_n)$ in $Z$ , and sequences of numbers
$(\tau_{n})\to \infty$ and $(\nu _{n})$ such that for all $n$
we have $d(\gamma(\tau_n),\xi_n(\nu _n))\leq \kappa$.  
\newline 
{\it Claim 1:} $\nu_n\to \infty$. 
\newline 
{\it Proof:} Suppose not. Then
there is a number $B$ such that for all $n$ $d(\xi_n(\nu_n),x)\leq
B$. By the triangle inequality 
$d(\gamma(\tau_n),x)\leq \kappa +B$, 
contradicting the fact that $\tau_n\to \infty$.  
\newline 
We write $\zeta_n=min\{\tau_n,\nu_n\}$. Then $(\zeta_n)\to\infty$.
\newline 
{\it Claim 2:} $d(\gamma(\zeta_n),\xi_n(\zeta_n))\leq
\kappa $.  
\newline 
{\it Proof:} Without loss of generality assume
$\tau_n\leq \nu_n$. Consider the geodesic triangle 
$[x, \gamma(\tau_n), \xi_n(\nu_n)]$. By Comparison Triangles we have 
\begin{equation*}
d(\gamma(\tau_n),\xi_n(\tau_n))\leq d(\gamma(\tau_n),\xi_n(\nu_n))\leq
\kappa. 
\end{equation*} 
Here, $\tau_n=\zeta_n$ so the Claim is proved.
\newline 
Since $Z$ is compact we may assume $\{\xi_n\}$ converges to
some $\xi\in Z$. Choose a sequence of numbers  $(\delta_n)\to \infty$
such that $d(\xi_n(\delta_n),\xi(\delta_n))\leq 1$. Then $\xi_n$ is
uniformly within $1$ of $\xi$ on $[0,\delta_n]$. If $\zeta_n\geq \delta_n$
then $d(\gamma(\delta_n),\xi(\delta_n))\leq \kappa +1$ and in that case we
write $\eta_n:=\delta_n$. On the other hand, if $\zeta_n<\delta_n$,
then $d(\gamma(\zeta_n),\xi(\zeta_n))\leq \kappa +1$ and we write
$\eta_n:=\zeta_n$. Either way, $\eta_n\to \infty$ so $\gamma$ and
$\xi$ are asymptotic, both starting at $x$, meaning that they are identical. This
contradicts the fact that $c\not\in Z$.  
\end{proof}

\begin{thm}\label{weakcut}  If $X$ does not have semistable fundamental
group at infinity then there are points $a, b, c\in \del X$ such that $c$
weakly separates $a$ from $b$ in $\del X$.

\end{thm}
\begin{proof} Let $p\in X$. By Lemma \ref{semi} we know that there
exists $\rho >0$ such that for any integer $n>\rho $ there exists ${\sigma}_{n}>n$ and
unit speed geodesic rays $\alpha _{n}, \beta _{n}:[0,\infty)\to X$
starting at $p$ with $d(\alpha _{n}({\sigma}_{n}), \beta _{n}({\sigma}_{n}))\leq
2\lambda +1$,
but the segment $[\alpha _{n}({\sigma}_{n}), \beta _{n}({\sigma}_{n})]$ cannot be
homotoped relative to $\alpha _{n}$ and $\beta _{n}$ out of 
$B(p, {\sigma}_{n}+\lambda +1)$ by a homotopy missing $B(p, \rho)$.  
\vskip 5pt 
Choose a compact set $K$ whose $G$-translates cover $X$. For each $n$ choose an
isometry $g_{n}$ with $g_{n}(\alpha _{n}({\sigma}_{n}))\in K$. Replacing $K$ by its
$2\lambda +2$-neighborhood we may assume $g_{n}(\beta _{n}({\sigma}_{n}))\in K$
too. 
Passing to a subsequence we may assume: 
\begin{enumerate}[(i)]
\item $\{g_{n}p\}$ converges to a point $c\in\del X$;
\item $\{g_{n}(\alpha _{n}({\sigma}_{n}))\}$ converges to a point $x\in K$ and  
$\{g_{n}(\beta _{n}({\sigma}_{n}))\}$ converges to a point $y\in K$;
\item $\{g_{n}\alpha _{n}\}$ converges to $\alpha$ and $\{g_{n}\beta
_{n}\}$ converges to $\beta$, where $\alpha $ and $\beta $ are geodesic
lines with $\alpha (-\infty)=\beta(-\infty)=c$. Here, convergence is uniform on compact
subsets.
\end{enumerate}
\vskip 5pt
Parametrize $\alpha , \beta :{\R}\to X$ with unit
speed so that $\alpha (0)=x$ and $\beta (0)=y$. Note that $d(x,y)\leq 2\lambda +1$. Let
$a=\alpha (\infty )$ and let $b=\beta (\infty)$.
\vskip 5pt
We show that $a\neq b$ and that $c$ weakly separates $a$ from $b$. Suppose there is a
continuum $Z\subseteq \del X$ containing $a$ and $b$ but not containing\footnote{The proof explains why 
the existence of such a $Z$ contradicts the hypothesis of non-semistability; if $a$ were
equal to $b$ then $Z:=\{a\}$ would be such a continuum; thus the proof shows 
$a\neq b$.} $c$. By Lemma \ref{compact}, since $c\notin Z$ we may choose
$m$ large enough that the cone ${\mathcal C}_{x}Z$  is disjoint from
$B(g_{m}p, \rho +2\lambda +2)$.  We may also assume that the segments
$[g_{m}(\alpha _{m}({\sigma}_{m})),  g_{m}(\alpha _{m}({\sigma}_{m}+2\lambda
+2))]$ and  
$[g_{m}(\beta _{m}({\sigma}_{m})),  g_{m}(\beta _{m}({\sigma}_{m}+2\lambda
+2))]$ lie within distance $1$ of $[\alpha (0), \alpha (2\lambda +2)]$ and $[\beta (0),
\beta (2\lambda +2)]$ respectively.
\vskip 5pt     

For $d\in Z$ let $\gamma _{d}:[0,\infty )\to X$ be the geodesic ray
(``cone line") starting at $x$ and ending at $d$. There is $\kappa $
with $\gamma _{d}(\tau )\notin B(g_{m}p, {\sigma }_{m}+2\lambda +2)$
for all $\tau \geq \kappa$ and all $d\in Z$. Let $\widehat Z=\{\gamma
_{d}(\kappa )\mid d\in Z\}$. This is a connected subset of $X$ since it
is a continuous image of $Z$.  \vskip 5pt

For small $\epsilon$, there is a piecewise geodesic path $\theta $
from $\gamma _{a}(\kappa )$ (i.e.  $\alpha (\kappa )$) to $\gamma
_{b}(\kappa )$, where $\theta $ lies uniformly within $\epsilon $ of
$\widehat Z$. Certainly, $\epsilon $ can be chosen so that $\theta $
misses $B(g_{m}p, \sigma_{m}+2\lambda +2)$.  The geodesic segment $\eta
$ joining $\gamma _{b}(\kappa )$ to $\beta (\kappa )$ has length is at
most $2\lambda +1$. For a suitably large choice of $\kappa$ this also
misses $B(g_{m}p, \sigma_{m}+2\lambda +2)$.  The combined path $\theta
\star \eta$ joins $\alpha (\kappa )$ to $\beta (\kappa )$ 

\vskip 5pt

The straight line unit-speed contraction of $\theta $ to the point $x$
misses $B(g_{m}p, \rho +2\lambda +2)$, and during that contraction the
point $\gamma _{b}(\kappa )$ traverses the segment $[\gamma _{b}(\kappa
), x]$.  That segment lies within $2\lambda +1$ of the segment $[\beta
(\kappa ), y]$.  Consider the geodesic quadrilateral $[\beta (\kappa ),
y]\cup [y,x]\cup [x,\gamma _{b}(\kappa )]\cup [\gamma _{b}(\kappa ),
\beta (\kappa )]$.  The geodesic segment from each point on  $\eta $ to
the corresponding point on $[x, y]$ defines a homotopy between $\eta $
and $[x, y]$ which stays within $2\lambda +1$ of $[x, \gamma _{b}(\kappa
)]$. This can be matched with the previously defined contraction.
The result is a homotopy of $\theta \star \eta$ to $[x,y]$ relative
to $\alpha$ and $\beta $, which misses $B(g_{m}p,\rho +1)$ (see II.2.2
of \cite{BrHa99}). In particular, the segment $[x,y]$ can be homotoped
to a path outside $B(g_{m}p, \sigma_{m} +2\lambda +2)$ but inside $B(g_{m}p,
\sigma_{m} +2\lambda +3)$ by a homotopy, relative to $\alpha $ and $\beta $,
which misses $B(g_{m}p, \rho+1)$. Since the segments 
$[g_{m}(\alpha _{m}({\sigma}_{m})),  g_{m}(\alpha
_{m}({\sigma}_{m}+2\lambda +2))]$ and
$[g_{m}(\beta _{m}({\sigma}_{m})),  g_{m}(\beta _{m}({\sigma}_{m}+2\lambda
+2))]$ 
lie within distance $1$ of $[\alpha (0), \alpha (2\lambda
+2)]$ and $[\beta (0), \beta (2\lambda +2)]$, there is a ``nearby''
homotopy moving the segment $[g_{m}(\alpha_{m}({\sigma}_{m})),
g_{m}(\beta _{m}({\sigma}_{m}))]$, relative to $g_{m}\alpha _{m}$ and
$g_{m}\beta _{m}$, to a path outside $B(g_{m}p, \sigma _{m} +2\lambda +2)$
but inside $B(g_{m}p, \sigma _{m} +2\lambda +3)$ while missing $B(g_{m}p,
\rho)$. Translating back by $g_m^{-1}$ gives the required contradiction.
\end{proof} 

The following addendum strengthens Theorem \ref{weakcut}:

\begin{addendum} There are neighborhoods $U$ of $a$ and $V$ of $b$
(in $\del X$) such that $c$ weakly separates $U$ from $V$.  
\end{addendum}

\begin{proof} A suitable neighborhood $U$ of $a$ is defined by all
geodesic rays starting at $x$ that stay close to $\alpha$ on the interval
$[0, 2\lambda +4]$; the neighborhood $V$ is defined similarly with respect
to $\beta $ and $y$. Let $\alpha'$ and $\beta'$ be such geodesic rays
ending in $\del X$ in the points $a'$ and $b'$ respectively. Suppose
there is a continuum $Z'\subseteq \del X$ containing $a'$ and $b'$ but
not $c$. Then the argument in the proof of Theorem \ref{weakcut} gives
an appropriate homotopy of the segment $[x,y]$, relative to $\alpha'$
and $\beta'$, to a path outside $B(g_{m}p, \sigma_{m} +2\lambda +2)$ but
inside $B(g_{m}p, \sigma _{m} +2\lambda +3)$. As in that theorem one then
obtains a homotopy moving the segment $[g_{m}(\alpha_{m}({\sigma}_{m})),
g_{m}(\beta _{m}({\sigma}_{m}))]$, relative to $g_{m}\alpha _{m}$ and
$g_{m}\beta _{m}$, to a path outside $B(g_{m}p, \sigma_{m} +2\lambda
+2)$ but inside $B(g_{m}p, \sigma _{m} +2\lambda +3)$, and this homotopy
misses $B(g_{m}p, \rho)$ Applying $g_m^{-1}$ as before gives the required
contradiction.  
\end{proof}

\section{Appendix: Proof of Proposition \ref{tech} }\label{appendix}  
Recall that the space $X$, and hence also the space $\widehat X$ (obtained by adding a
$Z$-set
to $X$), are AR's. Let $Q$ denote the Hilbert Cube, equipped with the metric it inherits
as a
subspace of Hilbert Space. Then $X\times Q$, with the usual product metric, is also a
proper
$CAT(0)$ space. The boundary of $X\times Q$ is again $\del X$. If
$\{U_{n}:=X-K_{n}\}$ of
Section \ref{stab} is a basic sequence of neighborhoods of infinity in $X$ then 
$\{U_{n}\times
Q\}$ plays the same role in $X\times Q$. 
\vskip 5pt     
By theorems of West \cite{We} and Chapman \cite{Ch}, $X\times Q$ is a $Q$-manifold
and
hence is homeomorphic to $P\times Q$ for some locally compact polyhedron $P$.
Clearly the
polyhedral structure gives a sequence $\{V_{n}\}$ of closed ANR neighborhoods of
infinity in
$P\times Q$ (which become compact ANR's when compactified by the boundary). So the
shape
theory of $\del X $ can be read off from the pro-homotopy of the sequence $\{V_{n}\}$.
The
image, under the homeomorphism, of the $V_{n}$'s in $X\times Q$ gives a sequence of
closed
ANR neighborhoods of infinity in $X\times Q$. This image is cofinal in $\{U_{n}\times
Q\}$.
So one has semistable fundamental group if and only if the other has. 
\vskip 5pt
In summary, semistability of $\{V_{n}\}$ is equivalent to $\del X$ having semistable
fundamental pro-group, and is also equivalent to semistability of the inverse sequence of
fundamental groups of the $U_n\times Q$'s, hence also to semistability of the inverse
sequence
of fundamental groups of the $U_n$'s. \hfill$\square$
\vskip 5pt
\begin{Rem} It might be objected that the use of the deepest theorems of
infinite-dimensional
topology to prove the ``almost obvious" Proposition \ref{tech} is overkill. Other proofs
are
possible. One could follow a method used in \cite{On05} where the $CAT(0)$ structure
is used
to construct a good open cover of $X$ whose nerve is proper homotopy equivalent to
$X$. One
would then proceed, much as in our proof, to find closed ANR neighborhoods of infinity
in that
nerve. Another approach would be to use more general shape theory as described in the
monograph \cite{MS} (see especially pages 18-19) to prove that the sequence
$\{U_{n}\}$ is an
``expansion" of $\del X$. Open subsets of $X$ are ANR's and hence are homotopy
equivalent to
CW complexes, so this is a  ``polyhedral expansion" in their sense. This would establish
that the
shape theory of $\del X$ can be read off from $\{U_{n}\}$. 
\end{Rem}

\def\cprime{$'$}
\providecommand{\bysame}{\leavevmode\hbox to3em{\hrulefill}\thinspace}
\providecommand{\MR}{\relax\ifhmode\unskip\space\fi MR }
\providecommand{\MRhref}[2]{%
  \href{http://www.ams.org/mathscinet-getitem?mr=#1}{#2}
}
\providecommand{\href}[2]{#2}


\begin{thebibliography}{MsS82}

\bibitem[BB08]{BaBu}
Werner Ballmann and Sergei Buyalo, \emph{Periodic rank one geodesics in
  {H}adamard spaces}, Geometric and probabilistic structures in dynamics,
  Contemp. Math., vol. 469, Amer. Math. Soc., Providence, RI, 2008, pp.~19--27.
  \MR{2478464}

\bibitem[BH99]{BrHa99}
Martin~R. Bridson and Andr{\'e} Haefliger, \emph{Metric spaces of non-positive
  curvature}, Grundlehren der Mathematischen Wissenschaften [Fundamental
  Principles of Mathematical Sciences], vol. 319, Springer-Verlag, Berlin,
  1999. \MR{1744486 (2000k:53038)}

\bibitem[Cha74]{Ch}
T.~A. Chapman, \emph{Topological invariance of {W}hitehead torsion}, Amer. J.
  Math. \textbf{96} (1974), 488--497. \MR{0391109}

\bibitem[Geo08]{G}
Ross Geoghegan, \emph{Topological methods in group theory}, Graduate Texts in
  Mathematics, vol. 243, Springer, New York, 2008. \MR{2365352}

\bibitem[GK91]{GK}
Ross Geoghegan and J{\'o}zef Krasinkiewicz, \emph{Empty components in strong
  shape theory}, Topology Appl. \textbf{41} (1991), no.~3, 213--233.
  \MR{1135099}

\bibitem[KM79]{KM}
J\'ozef Krasinkiewicz and Piotr Minc, \emph{Generalized paths and pointed
  {$1$}-movability}, Fund. Math. \textbf{104} (1979), no.~2, 141--153.
  \MR{551664}

\bibitem[Kra77]{K}
J\'ozef Krasinkiewicz, \emph{Local connectedness and pointed 1-movability},
  Bull. Acad. Polon. Sci. Sé Sci. Math. Astronom. Phys. \textbf{25} (1977),
  no.~12, 1265--1269. \MR{500986}

\bibitem[MsS82]{MS}
Sibe Marde\v~si\'c and Jack Segal, \emph{Shape theory}, North-Holland
  Mathematical Library, vol.~26, North-Holland Publishing Co., Amsterdam-New
  York, 1982, The inverse system approach. \MR{676973}

\bibitem[Ont05]{On05}
Pedro Ontaneda, \emph{Cocompact {CAT}(0) spaces are almost geodesically
  complete}, Topology \textbf{44} (2005), no.~1, 47--62. \MR{2104000
  (2005m:57002)}

\bibitem[PS09]{PS}
Panos Papasoglu and Eric Swenson, \emph{Boundaries and {JSJ} decompositions of
  {CAT}(0)-groups}, Geom. Funct. Anal. \textbf{19} (2009), no.~2, 559--590.
  \MR{2545250}

\bibitem[Wes77]{We}
James~E. West, \emph{Mapping {H}ilbert cube manifolds to {ANR}'s: a solution of
  a conjecture of {B}orsuk}, Ann. of Math. (2) \textbf{106} (1977), no.~1,
  1--18. \MR{0451247}

\end{thebibliography}
\end{document}